\documentclass[12pt]{article}
\usepackage{amsfonts}
\newtheorem{theorem}{Theorem}[section]
\newtheorem{lemma}[theorem]{Lemma}
\newtheorem{cor}[theorem]{Corollary}
\newtheorem{prop}[theorem]{Proposition}
\newtheorem{example}{Example}
\newenvironment{eg}{\begin{example}\normalfont}{\end{example}}

\newtheorem{proof}{Proof}

\newenvironment{pf}{\begin{proof}\normalfont}{\quad$\square$\end{proof}}

\newcommand{\mv}{\mathbf{v}}

\begin{document}

\title{Equitable partitions of Latin-square graphs}
\author{R. A. Bailey\thanks{School of Mathematics and Statistics,
University of St Andrews, St Andrews, Fife KY16 9SS, UK}
\thanks{The first two authors are grateful to  Shanghai Jiao Tong University for
funding, from NSFC (11671258) and STCSM (17690740800), 
a research visit where part of this work was done.} , 
Peter J. Cameron$^{*\dag}$, 
Alexander L. Gavrilyuk\thanks{Center for Math Research and Education, 
Pusan National University,
2, Busandaehak-ro 63beon-gil, Geumjeong-gu, Busan, 46241, Republic of Korea;
and
Krasovskii Institute of Mathematics and Mechanics,
S. Kovalevskaya st., 16, Yekaterinburg 620990, Russia}
\thanks{Research supported by
BK21plus Center for Math Research and Education at PNU.}
\\ 
and Sergey V. Goryainov\thanks{Shanghai Jiao Tong University, 800 Dongchuan Rd., 
Minhang District, Shanghai, China;
Krasovskii Institute of Mathematics and Mechanics, 
S. Kovalevskaya st.,~16, Yekaterinburg 620990, Russia; and
Chelyabinsk State University, Brat'ev Kashirinyh st.,~129, Chelyabinsk, 454021, Russia}
\thanks{Research supported by NSFC (11671258), STCSM (17690740800) and
RFBR (17-51-560008).}}
\date{January 2018}
\maketitle

\begin{abstract}
We study equitable partitions of Latin-square graphs, and give a complete
classification of those whose quotient matrix does not have an eigenvalue $-3$.


MSC: 05 E 30; 05 C 50; 05 B 15
\end{abstract}

\section{Introduction}

In the International Workshop on Bannai--Ito Theory in Hangzhou in November
2017, the fourth author spoke about his result with the third author
classifying the equitable partitions of the bilinear-forms graph
$\mathrm{Bil}_2(2\times d)$ \cite{g_hangzhou}. This graph can be regarded as
the Latin-square graph associated with the Cayley table of the additive group
of the $d$-dimensional vector space $W$ over the $2$-element field. An
equitable partition is associated with a matrix whose spectrum is contained in
that of the adjacency matrix of the graph; the result was a complete 
classification in the case where the eigenvalue $-3$ of the adjacency matrix
does not occur. The parts of such a partition must be unions of rows, columns,
or letters of the Cayley table, or subsquares corresponding to subspaces of
$W$ of codimension~$1$.

Here we present an extension of the result, to construct and classify all
equitable paritions of arbitrary Latin-square graphs for which the eigenvalue
$-3$ (the smallest of the three eigenvalues of any Latin-square graph) does not
occur (see Theorem~\ref{t:main} below). Remarkably, a relatively small
generalisation of the subspace construction is required; we replace this by
the notion of an inflation of a ``corner set'' in the Cayley table of a cyclic
group.

We begin in Section~\ref{s:gen}
with some preliminaries about equitable partitions, including showing that, for
the main theorem, it is enough to classify sets of vertices which are parts
of $2$-part equitable partitions, and are minimal subject to this condition.
We define Latin-square graphs in Section~\ref{s:LSG},
and in Section~\ref{s:blowup} we give the inflation construction which
is used to produce examples. 
Section~\ref{s:main} contains the main substance of this paper:
we construct the examples, and prove that there are no more.

The paper concludes with three short sections discussing possible further
directions. Sections~\ref{sec:B}--\ref{s:mixed} explore equitable partitions
which do involve the eigenvalue $-3$; we give examples to show that
classification in this case is unlikely to be feasible. 
If the Latin square has order $n$ then the other non-principal eigenvalue is $n-3$
(see Section~\ref{s:gen}):
we observe that examples which do not involve this eigenvalue
are connected with orthogonal arrays.
Finally, Section~\ref{s:mols} generalizes Latin squares to sets of mutually orthogonal Latin
squares.

Equitable partitions of distance-regular graphs involving only the two largest
eigenvalues of the graph were considered by Meyerowitz~\cite{meyerowitz}, who
classified them for Hamming and Johnson graphs. In general, the classification
of equitable partitions seems a hard problem, since it includes questions such
as tight sets in polar spaces, see \cite{bcd,bklp,blp}. Another classification
result is given in \cite{gg}.

\section{Equitable partitions}
\label{s:gen}

A partition $\Delta=\{\Delta_1,\ldots,\Delta_r\}$ of the vertex set of a
graph $\Gamma$ is said to be \emph{equitable} if there is an $r \times r$ 
matrix $M=(m_{ij})$ 
such that the number of vertices of $\Delta_j$
joined to a vertex $\omega\in\Delta_i$ is $m_{ij}$, depending on $i$ and $j$
but not on the choice of $\omega$. This term is used by Godsil and
Royle~\cite[\S9.3]{gr}; Fon-Der-Flaass~\cite{fdf} called such partitions
\emph{perfect}. We shall reserve this term for a set which is a part of a
$2$-part equitable partition, see below.

The spectrum of $M$ is contained in the spectrum of the adjacency matrix
$A(\Gamma)$ of the graph $\Gamma$: indeed, the characteristic polynomial of
$M$ divides that of $A(\Gamma)$~\cite[Theorem 9.3.3]{gr}. Since this result is
crucial to our approach, we outline a proof. The matrix $M$ is called the
\emph{quotient matrix} of the equitable partition. When we speak of eigenvalues
of an equitable partition, we refer to eigenvalues of the corresponding
quotient matrix.

We begin with some general information about equitable partitions. Let $\Omega$
be the vertex set of $\Gamma$, with $|\Omega|=N$;
and let $V$ be the $N$-dimensional vector space $\mathbb{R}^\Omega$, whose
basis vectors correspond to the vertices of $\Gamma$. Let $A(\Gamma)$ be the
adjacency matrix of $\Gamma$, and let $\mv_i\in V$ be the characteristic vector
of the part $\Delta_i$.

From the definition of an equitable partition, we see that
\[\mv_iA(\Gamma)=\sum_{j=1}^rm_{ji}\mv_j,\]
so that the space $W=\langle \mv_1,\ldots,\mv_r\rangle$ is invariant under 
$A(\Gamma)$, and the restriction of $A(\Gamma)$ to this subspace has matrix $M$ 
relative to the given basis. (Indeed this property is equivalent to the
partition being equitable.) Hence, if the (pairwise orthogonal) eigenspaces of
$A(\Gamma)$ are $V_1,\ldots,V_e$, then
\[W=(W\cap V_1)\oplus\cdots\oplus(W\cap V_e),\]
so the spectrum of $M$ is contained in that of $A(\Gamma)$, and the cited
result follows.

\medskip

From now on, we assume that $\Gamma$ is a connected regular graph with 
valency $k$. Then $k$ is a simple eigenvalue of $A(\Gamma)$. Moreover, the quotient
matrix~$M$
of an equitable partition has all row sums equal to $k$, so that $k$ is an 
eigenvalue of $M$. We call $k$ the \emph{principal eigenvalue}.

We say that an equitable partition $\Delta$ is \emph{$\mu$-equitable} if its
quotient matrix $M$ has all 
non-principal eigenvalues equal to $\mu$. Furthermore, we
call a  non-empty proper subset $S$ of $\Omega$ a \emph{$\mu$-perfect set} if
the partition $\{S,\Omega\setminus S\}$ is $\mu$-equitable.
Note that, if a set $S$ is
$\mu$-perfect, then so is its complement $\Omega\setminus S$.

\begin{prop}
Let $\Delta=\{\Delta_1,\ldots,\Delta_r\}$ be a partition of the vertex set
$\Omega$ of the regular connected graph $\Gamma$.
\begin{enumerate}\itemsep0pt
\item If $\Delta$ is $\mu$-equitable, then each set $\Delta_i$ is $\mu$-perfect.
\item Conversely, if $\Delta_1,\ldots,\Delta_{r-1}$ are all $\mu$-perfect,
then $\Delta$ is $\mu$-equitable.
\end{enumerate}
\label{p:equit}
\end{prop}

\begin{pf}
(a) Suppose that the hypotheses hold, and let $\mv_i$ be the characteristic
vector of $\Delta_i$. Then $\mv_i$ lies in the space $V_0\oplus V_1$, where
$V_0$ is the $k$-eigenspace (spanned by the all-$1$ vector $\mv_0$) and $V_1$
the $\mu$-eigenspace. Hence the span of $\mv_0$ and $\mv_i$ is
$A(\Gamma)$-invariant, and the restriction of $A(\Gamma)$ to this subspace has
eigenvalues $k$ and $\mu$; thus $\Delta_i$ is a $\mu$-perfect set.

\smallskip

(b) Conversely suppose that $\Delta_1,\ldots,\Delta_{r-1}$ are
$\mu$-perfect. Then the subspace spanned by $\mv_1,\ldots,\mv_{r-1}$ and
$\mv_0$ is contained in $V_0\oplus V_1$; since this space also contains
$\mv_r$, the conclusion follows.
\end{pf}

\begin{cor}
Let $S$ be a $\mu$-perfect set, and $T$ a non-empty proper subset of $\Omega\setminus S$.
Then $T$ is $\mu$-perfect if and only if $S\cup T$ is $\mu$-perfect.
\label{c:2}
\end{cor}

\begin{pf}
The forward direction follows immediately from 
Proposition~\ref{p:equit}(b):
if $S$ and $T$ are $\mu$-perfect then $\{S,T,\Omega\setminus(S\cup T)\}$
is $\mu$-equitable, 
and so
$\Omega\setminus(S\cup T)$ (and also its complement $S\cup T$) is
$\mu$-perfect. For the converse, if $S\cup T$ is $\mu$-perfect, then
two parts of the partition $\{S,T,\Omega\setminus(S\cup T)\}$ are
$\mu$-perfect; so this partition is $\mu$-equitable, and 
all its parts are $\mu$-perfect.
\end{pf}

\begin{cor}
If $\Delta$ is a $\mu$-equitable partition then any non-trivial coarsening of $\Delta$
is $\mu$-equitable.
\end{cor}

\begin{pf}
All parts of $\Delta$ are $\mu$-perfect, and so by 
Corollary~\ref{c:2}
the same is true for any non-trivial coarsening of $\Delta$; then
Proposition~\ref{p:equit}(b) applies.
\end{pf}

\section{Latin-square graphs}
\label{s:LSG}
Let $\Lambda$ be a Latin square of order $n$. Take $\Omega$ to be the set of
cells of $\Lambda$, so that $N = \left|\Omega\right |=n^2$. 
There are three uniform partitions
$R$, $C$ and $L$ of $\Omega$ into $n$ parts of size $n$. The parts of $R$ are
rows, the parts of $C$ are columns, and the parts of $L$ are letters.

If $\omega\in\Omega$, then $R(\omega)$, $C(\omega)$ and $L(\omega)$ denote
the row, column and letter containing $\omega$ (regarded as subsets of
$\Omega$).

As in Section~\ref{s:gen}, we denote by $V$
the $n^2$-dimensional vector space $\mathbb{R}^\Omega$.
Let $V_0$ be its one-dimensional subspace of constant vectors.  The
characteristic vectors of all rows span an $n$-dimensional subspace $V_R$
containing $V_0$.  Columns define a similar subspace $V_C$, and letters a
similar subspace $V_L$. Put $V_1 = (V_R + V_C + V_L) \cap  V_0^\perp$ and
$V_2 = (V_R + V_C + V_L)^\perp$, so that $V$ is the orthogonal direct sum of
$V_0$, $V_1$ and $V_2$.

The Latin square $\Lambda$ defines a Latin-square graph $\Gamma$ with
vertex set $\Omega$ and valency $k=3(n-1)$.  Each vertex is joined to every
other vertex in the same row or column or letter.  Denote the adjacency matrix
of $\Gamma$ by $A$.
We refer to the elements of $\Omega$ as cells or vertices, depending on the context.

The graph $\Gamma$ is strongly regular, so the matrices $I$, $A$ and $J-A-I$
form the adjacency matrices of an association scheme of rank three,
where $I$ is the identity matrix of order $n^2$ and $J$ is the $n^2 \times n^2$ matrix 
whose  entries are all equal to $1$.
The
common eigenspaces are $V_0$ (of dimension $1$), $V_1$ (of dimension $3(n-1)$)
and $V_2$ (of dimension $(n-1)(n-2)$). The eigenvalues of $A$ on these three
spaces are respectively $k=3(n-1)$, $n-3$, and $-3$.

In the special case $n=2$, the Latin square graph is the complete graph
$K_4$, and the eigenspace $V_2$ does not occur (the formula above gives its
dimension as zero).

\section{Inflation}
\label{s:blowup}

Here is a construction that we shall use several times.  

Let $\Lambda_0$ be a $t \times t$ Latin square on $\Omega_0$.  Replace each
occurrence of letter $i$ by an $s \times s$ Latin square on an alphabet
$\mathcal{A}_i$, where $\mathcal{A}_i \cap \mathcal{A}_j = \emptyset$ if
$i\ne j$, to obtain a Latin square $\Lambda_1$ of order $st$.  There is no
requirement for the $t$ Latin squares on alphabet $\mathcal{A}_i$ to be the
same, or even isomorphic.

This construction gives an \emph{orthogonal block structure} on a set
$\Omega_1$ of size $(st)^2$~\cite{bailey}.  The non-trivial partitions are
\begin{itemize}
\item
rows ($R$), columns ($C$), letters ($L$), each with $st$ parts of size $st$;
\item
fat rows ($\tilde{R}$), fat columns ($\tilde{C}$), fat letters ($\tilde{L}$),
each with $t$ parts of size $s^2t$, corresponding to the rows, columns and
letters of $\Lambda_0$;
\item
subsquares ($Q$), with $t^2$ parts of size $s^2$, where $Q$ is the infimum of
every pair of $\tilde{R}$, $\tilde{C}$ and $\tilde{L}$.
\end{itemize}
Like every orthogonal block structure, this defines an association scheme on
$\Omega_1$. The partition $Q$ is inherent in this, so it defines a quotient
scheme on the set of parts of $Q$.  This quotient scheme is precisely the
original Latin square $\Lambda_0$.

We call $\Lambda_1$ an \emph{$s$-fold inflation} of $\Lambda_0$.

\begin{theorem}
The partition $Q$ of $\Omega_1$ is equitable for the Latin-square graph 
$\Gamma_1$ defined by the Latin square $\Lambda_1$. 
\label{t:inflation_equitable}
\end{theorem}

\begin{pf}
A vertex in a part of $Q$ is joined to $3(s-1)$ further vertices in that part,
since the induced subgraph is a Latin-square graph from a square of order $s$.
It is joined to $s$ vertices in each part of $Q$ in the same fat row, fat column, or fat
letter, and to no vertex in any other part of $Q$. So the partition is equitable.
Its quotient matrix has the form $M=3(s-1)I+sA$, where $A$ is the adjacency matrix of
the Latin-square graph corresponding to $\Lambda_0$. The eigenvalues of $A$
are $3(t-1)$, $t-3$ and $-3$; so the eigenvalues of $M$ are 
$3(s-1)+3s(t-1)=3(st-1)$, $3(s-1)+s(t-3)=st-3$, and $3(s-1)+s(-3)=-3$,
the correct values for an equitable partition of a Latin-square graph with
$n=st$. Moreover, their multiplicities are those of the Latin-square graph
from $\Lambda_0$, namely $1$, $3(t-1)$, and $(t-1)(t-2)$.
\end{pf}

\begin{eg}
When $t=2$, there is a unique Latin square $\Lambda_0$ on an underlying set
$\Omega_0$ of size four.  The corresponding graph $\Gamma_0$ is complete, and
so all partitions of $\Omega_0$ are equitable for it.  Their $s$-fold
inflations give equitable partitions of some Latin squares of order $2s$. 
In this case, the 
multiplicities stated in the proof of Theorem~\ref{t:inflation_equitable}
show that all
non-principal eigenvalues of the partition are $2s-3$.

Indeed, if a Latin square of order $2s$ contains a subsquare of order $s$, 
then it necessarily arises as an inflation of the order-$2$ square. This
includes Cayley tables of groups of order $2s$ having subgroups of order~$s$.
\end{eg}

If $\Delta_0$ is a partition of $\Omega_0$ then the $s$-fold inflation gives a
partition $\tilde{\Delta}_0$ of $\Omega_1$ with the same number of parts.

\begin{theorem}
Let $\Lambda$ be an $s$-fold inflation of a Latin square $\Lambda_0$ of order $t$.
Let $\Gamma$ and $\Gamma_0$ be  the Latin-square graphs 
defined by $\Lambda$ and $\Lambda_0$ respectively.
\begin{enumerate}\itemsep0pt
\item
If $\Delta_0$ is
an equitable partition for $\Gamma_0$
then $\tilde{\Delta}_0$ is equitable for $\Gamma$.
\item 
If $P$ is a $(t-3)$-perfect subset of the 
vertex set of $\Gamma_0$
then the union of the $Q$-parts corresponding to the cells in $P$ is
an $(st-3)$-perfect subset of 
the vertex set of $\Gamma$.
\end{enumerate}
\label{thm:lift}
\end{theorem}

The proof is almost identical to that just given.

\section{$(n-3)$-perfect sets}
\label{s:main}

We return to the case where $\Gamma$ is the graph defined by a Latin square of
order~$n$.
Our goal is to describe the equitable partitions of $\Gamma$, especially those
with all non-principal eigenvalues equal to $n-3$. The preliminary results
we have given about equitable partitions show that every part of such a partition
is an $(n-3)$-perfect set, and any partition all of whose parts are
$(n-3)$-perfect is $(n-3)$-equitable.
So our job is to describe the $(n-3)$-perfect sets.

Suppose that an equitable $2$-partition has quotient matrix
\[
M = \left[\begin{array}{cc}p & b\\a & q\end{array}\right].
\]
Then 
\begin{equation}
p+b = a+q = k = 3(n-1).
\label{eq:1}
\end{equation}
Furthermore, if the non-principal eigenvalue is $n-3$ then
\begin{equation}
p+q = k+n-3 = 4n-6.
\label{eq:2}
\end{equation}

Moreover, if the first part is $S$, then counting edges between
$S$ and its complement gives $\left|S\right|b=(n^2-\left|S\right|)a$, 
so (since $a+b=2n$ from Equations~(\ref{eq:1}) and~(\ref{eq:2}))
we have
\begin{equation}
2\left|S\right|=na.
\label{eq:3}
\end{equation}

\subsection{Construction 1}

\begin{prop}
Any row, column or letter is an $(n-3)$-perfect set.
\end{prop}

\begin{pf}
Let $S$ be a row. (The other cases are similar.) The induced subgraph on $S$
is complete, and so any vertex in $S$ is joined to $2(n-1)$ vertices outside~$S$;
and any vertex outside $S$ is joined to two vertices of $S$ (one with the same
column and one with the same letter). So $\{S,\Omega\setminus S\}$ is
equitable, and its quotient matrix is
\[
M = \left[\begin{array}{cc}n-1 & 2(n-1)\\2 & 3n-5\end{array}\right].
\]
Thus the trace of $M$ is $4n-6$; since it has an eigenvalue $k=3n-3$, the
other eigenvalue is $n-3$, as required.
\end{pf}

It follows that any set which is a union of rows, or of columns, or of letters,
is $(n-3)$-perfect; and hence any partition all of whose parts are of this
form is equitable.

Another consequence of Corollary~\ref{c:2} is the following:

\begin{cor}
If an $(n-3)$-perfect set $S$ properly contains a row $T$, then $S\setminus T$ is
$(n-3)$-perfect; and similarly for a column or letter.
\end{cor}

So, in our search for the $(n-3)$-perfect sets, we may assume without loss that
such a set contains no row, column, or letter. We will call such a set
\emph{slender}.

\subsection{Slender sets and slices}

Given a slender subset $S$ of $\Omega$, call a \emph{slice} the
intersection of $S$ with any row, column or letter of $\Lambda$. 
Now we introduce some notation
for the size of a slice.  For each vertex $\omega$ in $\Omega$, put
$\rho(\omega) = \left|R(\omega) \cap S\right|$,
$\kappa(\omega) = \left|C(\omega) \cap S\right|$,
and
$\lambda(\omega) = \left|L(\omega) \cap S\right|$.
Then
\begin{equation}
\rho(\omega) + \kappa(\omega) + \lambda(\omega) = a 
\qquad \mbox{if $\omega\notin S$}.
\label{eq:5}
\end{equation}
In particular, no slice has size greater than $a$.  Also,
 Equations~(\ref{eq:1})--(\ref{eq:2}) show that
\begin{equation}
\rho(\omega) + \kappa(\omega) + \lambda(\omega) = 3 + p = k + n-q = n+a 
\qquad \mbox{if $\omega\in S$.}
\label{eq:6}
\end{equation}
Equation~(\ref{eq:6}) shows that if $\omega\in S$ then at least one of
$\rho(\omega)$,
$\kappa(\omega)$ and $\lambda(\omega)$ is greater than or equal to $(n+a)/3$.

\subsection{Construction 2}

Here is a construction of an equitable partition with three parts, two of which
are slender sets.

Let $\Lambda$ be the (back-)cyclic Latin square of order $n$, the Cayley
table of the cyclic group $Z_n$ of order $n$. We take the rows, columns,
and letters to be indexed by the set $\{0,1,\ldots,n-1\}$ of integers mod~$n$,
so that the letter in row~$i$ and column~$j$ is $i+j$ (with addition mod~$n$).

Consider the partition $\Delta$ with three parts $\Delta_{-1}$, $\Delta_0$
and $\Delta_1$, consisting of the cells $(i,j)$ with $i+j<n-1$, $i+j=n-1$, and
$i+j>n-1$ respectively (using integer addition here).

Figure~\ref{fig:back} shows the partition for $n=5$ and $n=6$, with 
$\Delta_{-1}$ in bold and $\Delta_0$ in calligraphic font.
For ease of reading, the letters indexed by $0$, $1$, $2$, etc.\ are shown as $A$, $B$, $C$,
etc.

\begin{figure}[htbp]
  \centering
  \begin{tabular}{c@{\qquad}c}
    $\begin{array}{|c|c|c|c|c|}
      \hline
      \mathbf{A} & \mathbf{B} & \mathbf{C} & \mathbf{D} & \mathcal{E}\\
      \hline
      \mathbf{B} & \mathbf{C} & \mathbf{D} & \mathcal{E} & A\\
      \hline
      \mathbf{C} & \mathbf{D} & \mathcal{E} & A & B\\
      \hline
      \mathbf{D} & \mathcal{E} & A & B & C\\
      \hline
      \mathcal{E} & A & B & C & D\\
      \hline
    \end{array}
    $ &
    $\begin{array}{|c|c|c|c|c|c|}
      \hline
      \mathbf{A} & \mathbf{B} & \mathbf{C} & \mathbf{D} & \mathbf{E} & \mathcal{F}\\
      \hline
      \mathbf{B} & \mathbf{C} & \mathbf{D} & \mathbf{E} & \mathcal{F} & A\\
      \hline
      \mathbf{C} & \mathbf{D} & \mathbf{E} & \mathcal{F} & A & B\\
      \hline
      \mathbf{D} & \mathbf{E} & \mathcal{F} & A & B & C\\
      \hline
      \mathbf{E} & \mathcal{F} & A & B & C & D\\
      \hline
      \mathcal{F} & A & B & C & D & E\\
      \hline
      \end{array}$
    \end{tabular}
  \caption{Equitable partitions of cyclic Latin squares, using Construction 2}
  \label{fig:back}
  \end{figure}

\begin{theorem}
With the above notation, the partition $\Delta$ is equitable, with both
non-principal eigenvalues equal to $n-3$.
\label{thm:cyc}
\end{theorem}

\begin{pf} We prove this by direct counting, see Figure~\ref{f:count}.

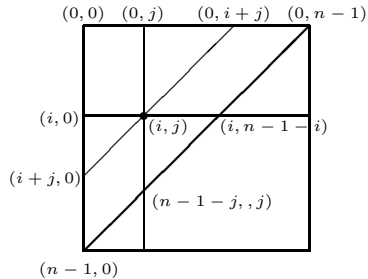
\begin{figure}[htbp]
\begin{center}
\setlength{\unitlength}{1mm}
\begin{picture}(30,30)
\thicklines
\multiput(0,0)(30,0){2}{\line(0,1){30}}
\multiput(0,0)(0,30){2}{\line(1,0){30}}
\put(0,0){\line(1,1){30}}
\thinlines
\put(8,18){\circle*{1}}
\put(8,0){\line(0,1){30}}
\put(0,18){\line(1,0){30}}
\put(0,10){\line(1,1){20}}
\put(-3,31){\tiny$(0,0)$}
\put(5,31){\tiny$(0,j)$}
\put(15,31){\tiny$(0,i+j)$}
\put(27,31){\tiny$(0,n-1)$}
\put(-6,17){\tiny$(i,0)$}
\put(8.5,16){\tiny$(i,j)$}
\put(18,16){\tiny$(i,n-1-i)$}
\put(-10,9){\tiny$(i+j,0)$}
\put(9,6){\tiny$(n-1-j,,j)$}
\put(-6,-3){\tiny$(n-1,0)$}
\end{picture}
\end{center}
\caption{\label{f:count}Counting neighbours in Theorem~\ref{thm:cyc}}
\end{figure}

Take a cell $(i,j)$ in $\Delta_{-1}$. Within $\Delta_{-1}$, there are $n-i-2$
cells in the same row, $n-j-2$ in the same column, and $i+j$ with the same
letter (excluding the cell $(i,j)$ itself); so it has $2n-4$ neighbours in
$\Delta_{-1}$. The cells in $\Delta_0$ all have letter $n-1$, which never occurs
in $\Delta_{-1}$; so $(i,j)$ is joined to two cells in $\Delta_0$, namely
$(i,n-1-i)$ (in the same row) and $(n-1-j,j)$ (in the same column). The 
remaining $n-1$ neighbours are in $\Delta_1$.

A cell $(i,n-1-i)$ in $\Delta_0$ is joined to the other $n-1$ cells in
$\Delta_0$ (all have the same letter), and to $n-1$ cells in $\Delta_{-1}$
(of which $n-i-1$ are in the same row and $i$ in the same column).

The other matrix coefficients follow by symmetry between $\Delta_{-1}$ and $\Delta_1$.

Thus the partition is equitable, with quotient matrix
\[M=\left[\begin{array}{ccc}
2n-4 & 2 & n-1\\ n-1 & n-1 & n-1\\ n-1 & 2 & 2n-4
\end{array}\right].\]
This matrix has trace $5n-9$; so its eigenvalues are $3n-3$, $n-3$, $n-3$.
\end{pf}

In particular, the parts of the partition are $(n-3)$-perfect. The part
$\Delta_0$ is a letter, but the other two parts are obviously slender.
Further partitions of this type can be found by changing the roles of rows,
columns, and letters.

We will call sets of the form $\Delta_{-1}$ in this example, possibly after
re-labelling of rows, columns and letters, \emph{corner sets}. Note that a
corner set is disjoint from a row and a column as well as a letter, and so
is a part of three different partitions of this type: as well as the one
given, we have
$\{\Delta_{-1}, \Delta_2, (\Delta_0\cup\Delta_1)\setminus\Delta_2\}$,
where $\Delta_2$ is either the last row or the last column.

Theorem~\ref{thm:lift} now shows that if $\Lambda$ is an $s$-fold inflation
of a cyclic Latin square of order $t$ then the inflation of the partition
$\Delta$ in Theorem~\ref{thm:cyc} is equitable with non-principal
eigenvalues $n-3$, where $n$ is the order of $\Lambda$.

In particular, a single cell in a Latin square of order $2$ is a corner set,
and inflation gives subsquares of order $n/2$ in Latin squares of even
order $n$. Hence such subsquares are $(n-3)$-perfect sets.

\subsection{The main theorem}

\begin{theorem}
Let $\Gamma$ be the Latin-square graph defined by a Latin square of order~$n$,
and $\Delta$ a partition of the vertex set of $\Gamma$. Then $\Delta$ is
$(n-3)$-equitable if and only if each part of $\Delta$ is a disjoint union of
rows, columns, letters, or inflations of corner sets.
\label{t:main}
\end{theorem}

We know from the results of Section~\ref{s:gen} that it is enough to
describe the $(n-3)$-perfect sets, and moreover that it is enough to show
the following:

\begin{theorem}
A slender $(n-3)$-perfect set in the Latin-square graph 
defined by a Latin square
of order $n$ is an inflation of a corner set.
\end{theorem}

The proof of the theorem is somewhat involved, so we begin with a 
summary and some comments on the notation.
Let $S$ denote a slender $(n-3)$-perfect set.

An inflation of a corner set, after suitable row and column permutations,
resembles the starred region shown in Figure~\ref{f:sketch}, which also shows fat rows
and fat columns.

\begin{figure}[htbp]
\begin{center}
\setlength{\unitlength}{1mm}
\begin{picture}(25,25)
\thinlines
\multiput(0,0)(0,5){6}{\line(1,0){25}}
\multiput(0,0)(5,0){6}{\line(0,1){25}}
\thicklines
\put(0,25){\line(1,0){20}}
\put(0,25){\line(0,-1){20}}
\multiput(0,5)(5,5){4}{\line(1,0){5}}
\multiput(5,5)(5,5){4}{\line(0,1){5}}
\multiput(1.5,21.5)(5,0){4}{$*$}
\multiput(1.5,16.5)(5,0){3}{$*$}
\multiput(1.5,11.5)(5,0){2}{$*$}
\put(1.5,6.5){$*$}
\end{picture}
\end{center}
\caption{\label{f:sketch}Sketch of the proof}
\end{figure}

We begin by identifying the subsquare in the top right of the figure: it is the intersection
of the fat row consisting of those rows meeting $S$ in the maximum number of columns and 
the fat column consisting of those columns not meeting $S$. 
Then, inductively, we work down and to the left, identifying the subsquares
on the boundary of $S$. 
The notation is introduced as the proof proceeds.

Parts (a)--(b) of the induction successively find rows whose slice sizes are strictly decreasing.
These rows are numbered by their position of occurrence in the
inductive proof rather than by their position in the square, either originally or after
permutation to the form in Figure~\ref{f:sketch}.  Row $i$ is the first row to be named
after fat row $\tilde{R}_{i-1}$ is identified in part (f), but it probably does not lie adjacent to
row $i-1$.

Part (c) of the induction uses row $i$ and fat column $\tilde{C}_1$ (on the right of 
Figure~\ref{f:sketch}) to define the fat letter $\tilde{L}_i$ and find some properties of it.
Then part~(d) uses row $i$ and fat letter $\tilde{L}_1$ to define fat column $\tilde{C}_i$ 
in such a way that $\tilde{L}_1$ is on the back-diagonal of the square.
Thus fat columns are numbered from right to left.  We do not really have a viable way of
numbering them that matches the reader's expectations, because we do not know at the start
of the proof that, for example, the size of the subsquares divides $n$.

The rest of parts (d) and (e) identify the letters in the intersection  of row $i$ with fat 
column $\tilde{C}_j$ for $j<i$.

Finally, part (f) shows that there is a fat row $\tilde{R}_i$ which contains row $i$ 
and which has the properties necessary for an inflated square.

Once the fat letters have been assigned to the subsquares outside $S$, 
the Latin square property forces their allocation 
to the subsquares in $S$, working upwards from the penultimate fat row.
However, the information gathered during the proof gives a more direct way of doing
this, as we show at the end of this section.

\medskip

We now embark on the details. 
Recall that $a$ is the (constant) number of neighbours in $S$ of any vertex outside $S$,
so that Equations~(\ref{eq:1})--(\ref{eq:6}) hold.
The proof makes frequent use of Equations~(\ref{eq:5}) and (\ref{eq:6}).

\begin{lemma}
\label{lem:maxslice}
If there is a slice of size $r$ then either $r=a$ or there is a slice of size
at least $(n+r)/2$.
\end{lemma}

\begin{pf}
Without loss of generality, assume that the slice of size $r$ is contained in
a row. Because $S$ is slender, there is a vertex $\alpha$ in this row
which is not in $S$.  Then Equation~(\ref{eq:5}) shows that
$\kappa(\alpha) + \lambda(\alpha) = a-r$.
If $r\ne a$ then at least one of $C(\alpha) \cap S$  and
$L(\alpha) \cap S$  is not empty. Without loss of generality, there is
a vertex $\beta$ in $C(\alpha) \cap S$ .
Then Equation~(\ref{eq:6}) shows that
$\rho(\beta) + \lambda(\beta) = n+a-\kappa(\beta) \geq (n+a) - (a-r) =n+r$.
Hence at least one of $\rho(\beta)$ and $\lambda(\beta)$ is at least $(n+r)/2$.
\end{pf}

\begin{cor}
There is at least one slice of size $a$.
\end{cor}

From this corollary and the fact that $S$ is slender, we see that $a<n$. We
put $s=n-a$.

\begin{figure}[htbp]
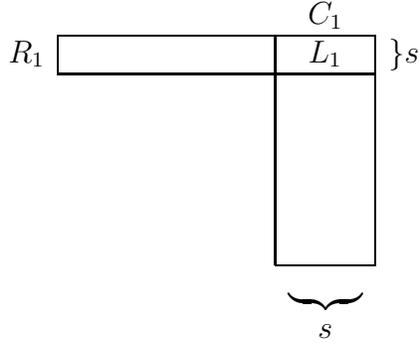

\[\begin{array}{c|c|c|c}
\multicolumn{2}{c}{} & \multicolumn{1}{c}{C_1} \\
\cline{2-3}
R_1 & \hspace{1in} & L_1 & \} s\\
\cline{2-3}
\multicolumn{2}{c|}{} &  &\\
\multicolumn{2}{c|}{} &  &\\
\multicolumn{2}{c|}{} &  &\\
\multicolumn{2}{c|}{} &  &\\
\multicolumn{2}{c|}{} &  &\\
\cline{3-3}
\multicolumn{2}{c}{} & \multicolumn{1}{c}{\underbrace{\hphantom{n-a}}}\\
\multicolumn{2}{c}{} & \multicolumn{1}{c}{s}
\end{array}
\]
\caption{What we know when $i=1$: all vertices in $\tilde{R}_1$ to the left
of the subsquare marked $L_1$ are in $S$, while all those in $\tilde{C}_1$ are outside $S$}
\label{fig:e1}
\end{figure}

Slightly abusing notation, write $\rho_i$, $\kappa_j$ and $\lambda_\ell$
for the size of the slice in row $i$, column $j$ and letter $\ell$
respectively.  Without loss of generality, we may assume that $\rho_1=a$.
Let $L_1$ be the set of $s$ letters whose cells in row~$1$ are not in $S$; and
let $C_1$ be the set of the $s$ columns whose intersection with
row~$1$ is not in $S$. 
The proof of Lemma~\ref{lem:maxslice} shows that  $\lambda(\alpha)=0=\kappa(\alpha)$ 
if $\alpha$ is in row~$1$ and a column in $C_1$. (See Figure~\ref{fig:e1}, which
also incorporates part~(f) of the following theorem for the case $i=1$.) 

The statement of the following theorem introduces further notation like $L_1$
and $C_1$ to denote various sets of letters, columns or rows.  In each case,
we use the same notation with a $\tilde{\hphantom{L}}$ on top to denote the
subset of $\Omega$ formed by the union of all letters or columns or rows in
that set.

\begin{theorem}
\label{thm:max}
Assume that $S$ is slender and $\rho_1=a$.
  If $n>(t-1)s$ and $1\leq i \leq t$ then the following are true.
Hence $n\geq ts$.  
\begin{enumerate}
  \item[(a)]
    There is no row $i$ with $n-is <\rho_i < n-(i-1)s$.
  \item[(b)]
    There is a row $i$ with $\rho_{i}= n-is$.
  \item[(c)]
    If $i>1$,  let $L_i$ be the set of $s$ letters in the intersection of
     row $i$    with $\tilde{C}_1$.
     If $i\geq 1$ then every letter $\ell$ in $L_i$ has $\lambda_\ell = (i-1)s$.
Hence $L_i$ and $L_{i'}$ are disjoint if $1\leq i'<i$.
\item[(d)]
If $i>1$,      let $C_{i}$ be the set of $s$ columns where row $i$ contains
    letters in $L_1$.  If $1\leq i'\leq i$ then
every vertex in the intersection of $\tilde{C}_{i'}$ with row $i$ is outside $S$.
  Hence, for every column $j$ in $C_{i}$, $\kappa_j=(i-1)s$.
  Moreover, every vertex in row $i$ outside
$\tilde{C}_1 \cup \cdots \cup \tilde{C}_i$  is in $S$.
  \item[(e)]
    If $1\leq j< i$ then the letters in the intersection of row $i$
    with $\tilde{C}_{j}$ are precisely those in $L_{i-j+1}$.
  \item[(f)]
    There are precisely $s$ rows whose slice has size $n- is$.
    If $R_i$ denotes the set of these rows, then $\tilde{R}_i \cap \tilde{C}_1$
    is a Latin square on the letters in $L_i$.  Moreover, if $1< j\leq i$ then
    $\tilde{R}_i\cap \tilde{C}_{j}$ is a Latin square on the letters in
    $L_{i-j+1}$. Also
  $\tilde{R}_i \setminus (\tilde{C}_1 \cup \cdots \cup \tilde{C}_i) \subset S$.
  \end{enumerate}
\end{theorem}

\begin{figure}[htbp]
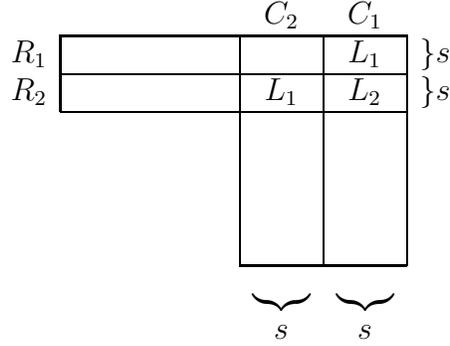

\[\begin{array}{c|c|c|c|c}
\multicolumn{2}{c}{} & \multicolumn{1}{c}{C_{2}} & \multicolumn{1}{c}{C_1} \\
\cline{2-4}
R_1 & \hspace{0.8in} & & L_1 & \} s\\
\cline{2-4}
R_2 & \hspace{0.8in} & L_1 & L_2 & \}s\\
\cline{2-4}
\multicolumn{2}{c|}{} &  &\\
\multicolumn{2}{c|}{} &  &\\
\multicolumn{2}{c|}{} &  &\\
\multicolumn{2}{c|}{} &  &\\
\cline{3-4}
\multicolumn{2}{c}{} & \multicolumn{1}{c}{\underbrace{\hphantom{s}}}
& \multicolumn{1}{c}{\underbrace{\hphantom{s}}}\\
\multicolumn{2}{c}{} & \multicolumn{1}{c}{s} & \multicolumn{1}{c}{s}
\end{array}
\]
\caption{What we know when $i=2$: all vertices to the left
of, or above,  subsquares marked $L_1$ are in $S$, while all those in, or below,
those subsquares are outside~$S$}
\label{fig:e2}
\end{figure}

\begin{pf}
  We use induction on $i$.  Parts (b), (d)  and (f) for $i=1$ give the
  situation summarized in
  Figure~\ref{fig:e1}. If $t\geq 2$, then parts (b) to (f) give
  Figure~\ref{fig:e2}
for $i=2$.

\begin{enumerate}
\item[(a)]
If $i=1$ this is true because no slice has size greater than $a$.
  
If $i>1$, assume that (a)--(f) are true for all $i'$ with $1\leq i'<i$.
These give the situation summarized in Figure~\ref{fig:e3}.

Let $\alpha$ be a vertex in 
$\tilde{C}_1 \setminus (\tilde{R}_1 \cup \cdots \cup \tilde{R}_{i-1})$.
Then $L(\alpha)$ is not in $L_1 \cup \cdots \cup L_{i-1}$, and so it
occurs in every row of $\tilde{R}_{i'}\cap S$ for $1\leq i' \leq i-1$, by (d) and (e).
By (f), $\lambda(\alpha) \geq (i-1)s$. Moreover, $\alpha\notin S$, and so
 $\rho(\alpha)+\lambda(\alpha)=a$.
Therefore $\rho(\alpha) \leq a- (i-1)s = n-is$.
This proves (a) for $i$.

\begin{figure}[htbp]
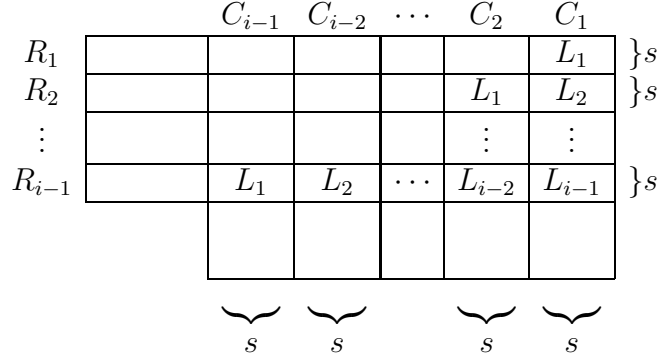

\[\begin{array}{c|c|c|c|c|c|c|c}
\multicolumn{2}{c}{} & \multicolumn{1}{c}{C_{i-1}} & \multicolumn{1}{c}{C_{i-2}} 
& \multicolumn{1}{c}{\cdots} &
\multicolumn{1}{c}{C_{2}} & \multicolumn{1}{c}{C_1} \\
\cline{2-7}
R_1 & \hspace{0.5in} & & & & & L_1 & \} s\\
\cline{2-7}
R_2 & \hspace{0.5in} & & & & L_1 & L_2 & \}s\\
\cline{2-7}
\vdots & & &  &&\vdots & \vdots&\\
\cline{2-7}
R_{i-1}  &\hspace{0.5in} & L_1 & L_2 & \cdots & L_{i-2} & L_{i-1} & \} s\\
\cline{2-7}
\multicolumn{2}{c|}{} &  & & & & \\
\multicolumn{2}{c|}{} &  & & & &\\
\cline{3-7}
\multicolumn{2}{c}{} & \multicolumn{1}{c}{\underbrace{\hphantom{s}}}
& \multicolumn{1}{c}{\underbrace{\hphantom{s}}}
& \multicolumn{1}{c}{} & \multicolumn{1}{c}{\underbrace{\hphantom{s}}}
& \multicolumn{1}{c}{\underbrace{\hphantom{s}}}\\
\multicolumn{2}{c}{} & \multicolumn{1}{c}{s} & \multicolumn{1}{c}{s} &
\multicolumn{1}{c}{} & \multicolumn{1}{c}{s} & \multicolumn{1}{c}{s} 
\end{array}
\]
\caption{What we know after $i-1$ steps in the induction: all vertices  to the
  left of, or above,  subsquares marked $L_1$ are in $S$,  while all those in,
  or below,
those subsquares are outside $S$}
\label{fig:e3}
\end{figure}

\item[(b)]
Since $\rho_1= a = n-s$, this is true when $i=1$.

If $i>1$, assume that (a)--(f) are true for all $i'$ with $1\leq i'<i$.
Suppose that the vertex $\alpha$ in the proof of~(a) is chosen to
minimize $\lambda(\alpha)$.  If $\lambda(\alpha) > (i-1)s$ then
$\rho(\alpha) < a-(i-1)s$, and so there is a vertex $\beta$ in
$R(\alpha) \setminus S \setminus(\tilde{L}_1 \cup \cdots \cup \tilde{L}_{i-1}) 
\setminus \tilde{C}_1$.
Then $\kappa(\beta)>0$ and
$a= \rho(\beta) + \kappa(\beta) + \lambda(\beta) =
\rho(\alpha) + \kappa(\beta) + \lambda(\beta) =
a-\lambda(\alpha)
+ \kappa(\beta) + \lambda(\beta)$.
It follows that $\lambda(\alpha)>\lambda(\beta)$.  But $\beta$ is not in
$\tilde{L}_1 \cup \cdots \cup \tilde{L}_{i-1}$, so there is some vertex $\gamma$
in $\tilde{C}_1 \setminus (\tilde{R}_1 \cup \cdots \cup \tilde{R}_{i-1})$ with 
$L(\gamma) = L(\beta)$.
This contradicts the choice of $\alpha$ to minimize $\lambda(\alpha)$.
It follows that $\lambda(\alpha)=(i-1)s$ and so $\rho(\alpha)=
a-(i-1)s = n - is$.  This proves (b) for $i$.

\item[(c)]
Assume that (b) is true for $i$.
  If $\ell\in L_i$ then there is a vertex $\alpha$ in $\tilde{C}_1$ with $\rho(\alpha)
  = n-is$ and $L(\alpha)=\ell$.  Since $\kappa(\alpha)=0$ and 
$\rho(\alpha)+\kappa(\alpha)   +\lambda(\alpha)=a$, this shows that 
$\lambda_\ell =\lambda(\alpha)= (i-1)s$.  

If $i>1$, assume that (c) is true for all $i'$ with $1\leq i'<i$.  
Then, for any such $i'$, any letter $m$ in $L_{i'}$ has $\lambda_m = (i'-1)s$.
Therefore $L_i$ and $L_{i'}$ are disjoint.
Thus (c) is true for $i$.

\item[(d)]
  Assume that (b) and (c) are true for $i$, and that (d) is true for all
  $i'$ with $1\leq i' <i$.
  If $1\leq i'<i$ and $\alpha$ is in $\tilde{C}_{i'}$ and row $i$
  then $\rho(\alpha)=n-is$ and $\kappa(\alpha)=(i'-1)s$.
  Since $i'<i$, we have $\rho(\alpha)+\kappa(\alpha) \leq n-2s<n$.
  If $\alpha\in S$ then   $\rho(\alpha)+\kappa(\alpha)
  +\lambda(\alpha)=n+a$, and so $\lambda(\alpha) >a$.
  This cannot happen,  and so $\alpha\notin S$.
  
  If $\alpha$ is in $\tilde{C}_{i}$ and row $i$ then $\alpha$ is in $\tilde{L}_1$.
Then   $\rho(\alpha) = n-is$ and $\lambda(\alpha)=0$. 
If $\alpha\in S$ then   $\rho(\alpha)+\kappa(\alpha)
  +\lambda(\alpha)=n+a$, and so $\kappa(\alpha) = a+is>a$.
This cannot happen, and so  $\alpha\notin S$.
Therefore
  $\rho(\alpha)+\kappa(\alpha)  +\lambda(\alpha)=a$, which
shows that $\kappa(\alpha)= (i-1)s$.

Finally, since  $\rho_{i}=n-is$ and all vertices in the intersection of
row $i$ and $\tilde{C}_1 \cup \cdots \cup  \tilde{C}_i$ are outside $S$, all the remaining
vertices in row $i$ must be in $S$.
Thus (d) is true for $i$.

\item[(e)]
If $i>1$, assume that (c), (d) and (f) are true for all $j$ with $1\leq j<i$.
If $\ell$ is a letter outside $L_1 \cup \cdots \cup L_{i-1}$ then
(d) shows that it occurs in $S$ in every row in $\tilde{R}_1 \cup \cdots \cup \tilde{R}_{i-1}$,
and so $\lambda_\ell \geq (i-1)s$.

If $1\leq j<i$ and $\alpha$ is a vertex in the intersection of row $i$ with
$\tilde{C}_{j}$, then (d) shows that $\alpha\notin S$.
Therefore $a = \rho(\alpha) + \kappa(\alpha) + \lambda(\alpha)
=  (n-is) + (j-1)s + \lambda(\alpha)$ and so $\lambda(\alpha) = (i-j)s$.
If $j=1$ then $i-j+1 = i$ and by definition the letters in the intersection of
row $i$ with $\tilde{C}_1$ are those in $L_{i}$. If $j>1$ then $\lambda(\alpha)
\leq (i-2)s$ and so $L(\alpha) \in L_1 \cup \cdots \cup L_{i-1}$.
Then it follows from (c) for integers less than
$i$ that $L(\alpha)\in L_{i-j+1}$.
Hence (e) is true for $i$.

\item[(f)]
Assume that (d) and (e) are true for $i$.

Let $\ell$ be a letter in $L_i$.  This occurs in $s$ rows of $\tilde{C}_1$,
all of whose slices have size $n-is$.
When $i=1$, each letter $m$ outside $L_1$ has
$\lambda_m >0$ and the argument in (c) shows that $m$ cannot occur
in the intersection of any of these rows with $\tilde{C}_1$.
For $i>1$, parts (d) and  (e) show that, for each of these rows,
the letters outside $S\cup \tilde{C}_1$ are precisely those in $L_1 \cup \cdots \cup L_{i-1}$.
Suppose that a letter $m$ in $L_i$ occurs on a vertex $\alpha$
in $S$ in such a row.
Then $\rho(\alpha)+\lambda(\alpha) = n-is + (i-1)s= n-s$
so $\kappa(\alpha) = n+a-n+s=n$, which is impossible because $S$
is slender.
Hence each of these $s$ rows intersects $\tilde{C}_1$ in 
a set of vertices whose letters
are the set $L_i$.

If there are any more rows with slice size $n-is$ then they must
contain each letter of $L_i$ in their slice.  The foregoing argument shows
that this cannot happen.

Now applying the arguments in (d) and (e) to each row in $R_i$
completes the proof of (f) for $i$.
\end{enumerate}
\end{pf}

Theorem~\ref{thm:max} shows that a Latin square $\Lambda$ with a slender
set has fat rows $\tilde{R}_i$, fat columns $\tilde{C}_j$ and fat letters
$\tilde{L}_\ell$, all of size $s$,
so that $n$ must be some multiple $ts$ of $s$.
Parts (c), (d) and (e) of the theorem explicitly assign $L_1$ to each 
intersection $\tilde{R}_i \cap \tilde{C}_i$ on the back-diagonal, and
$L_{i-j+1}$
to each intersection $\tilde{R}_i\cap \tilde{C}_j$  below the
back-diagonal.
If $\tilde{R}_i\cap \tilde{C}_j$ is above the back-diagonal  then 
$i<j$ because of the non-standard labelling of the fat columns.  Parts (b) and (d) show that
$\rho_i=n-is$ and $\kappa_j=(j-1)s$, so Equation~(\ref{eq:6}) gives
$\lambda(\alpha) = (t-j+i)s$ if $\alpha\in \tilde{R}_i\cap \tilde{C}_j$.
Then part (c) shows that  the letters  which occur in $\tilde{R}_i\cap \tilde{C}_j$ are precisely 
those in $L_{t-j+i+1}$.
Relabelling fat row $i$ as $i-1$, fat column $j$ as $t-j$ and fat letter $\ell$ as $\ell-2 \pmod t$
gives the back-cyclic Latin square of order $t$ in Construction~2.
Therefore $\Lambda$ is an $s$-fold inflation
of a back-cyclic Latin square of order $t$.

The elementary abelian $2$-group has no cyclic quotient of order greater than
$2$, and so the only inflation of a corner set which occurs in its Cayley
table is a subsquare corresponding to a subgroup of index~$2$. Thus we recover
the result of Gavrilyuk and Goryainov \cite{g_hangzhou} which was the starting point.

\section{$-3$-perfect sets}

\label{sec:B}
Let $\Delta_1$ be a non-empty proper subset of the set $\Omega$ of vertices of
a Latin-square graph, where $\left|\Omega\right|=n^2$ and
$\left|\Delta_1\right| = m$. Let $\Delta_2$ be the complement of $\Delta_1$,
so that $\left|\Delta_2\right| = n^2-m$. The \emph{contrast} between $\Delta_1$
and $\Delta_2$ is defined to be any non-zero multiple of the vector $\mathbf{z}$
which takes the value $n^2-m$ on each element of $\Delta_1$ and the value $-m$
on each element of $\Delta_2$. Now $\Delta_1$ and $\Delta_2$ are $-3$-perfect
sets if and only if this contrast is in $V_2$, which happens if and only if
the entries in $\mathbf{z}$ sum to zero on each row, column and letter.  This
means that the partition $\Delta$ is \emph{strictly orthogonal} to each of $R$,
$C$ and $L$: see~\cite[p.~8]{bailey_bcc}.
In the special case that $\left|\Delta_1\right| = \left|\Delta_2\right|$, this means that
$\{R,C,L,\Delta\}$ forms an \emph{orthogonal array of strength two} on $\Omega$.

If $\Delta_1$ is any \emph{transversal} for $\Lambda$ (a set of cells meeting
each row, column and letter just once) then $\Delta$ satisfies this condition.
More generally, $\Delta$ satisfies this condition if $\Delta_1$ is the union
of any collection of mutually disjoint transversals.  In particular, if there
is a Latin square $\Lambda'$ orthogonal to $\Lambda$ then any partition of its
letters gives a $-3$-equitable partition for $\Gamma$.

\begin{eg}
Figure~\ref{fig:4gl}(a) shows a Graeco-Latin square of order $4$.  Let
$\Lambda$ be the Latin square defined by the Latin letters. Let $\Delta_1$ be
the union of Greek letters $\alpha$ and $\beta$.  Then $\Delta_1$ and its
complement give the equitable partition for $\Gamma$ shown in
Figure~\ref{fig:4gl}(b).

\begin{figure}[htbp]
\centering
\begin{tabular}{c@{\qquad}c}
$\begin{array}{|cc|cc|cc|cc|}
\hline
A & \alpha & B & \beta & C & \gamma & D & \delta\\
\hline
B & \delta & A & \gamma & D & \beta & C & \alpha\\
\hline
C & \beta & D & \alpha & A & \delta & B & \gamma\\
\hline
D & \gamma & C & \delta & B & \alpha & A & \beta\\
\hline
\end{array}$
&
$\begin{array}{|c|c|c|c|}
\hline
\mathbf{A} & \mathbf{B} & C & D\\
\hline
B & A & \mathbf{D} & \mathbf{C}\\
\hline
\mathbf{C} & \mathbf{D} & A & B\\
\hline
D & C & \mathbf{B} & \mathbf{A}\\
\hline
\end{array}
$\\
\\
(a) & (b)
\end{tabular}
\caption{The Greek letters $\alpha$ and $\beta$ in the  Graeco-Latin square of order $4$ in
(a) give the equitable partition shown in (b), where the elements of $\Delta_1$ are shown in
bold}
\label{fig:4gl}
\end{figure}
\end{eg}

There seems to be no possibility of determining all transversals in Latin
squares. Even Ryser's celebrated conjecture, that any Latin square of odd order
contains a transversal, is still open. So there is no possibility for a
classification in this case similar to what we did for the eigenvalue $n-3$
in Section~\ref{s:main}.

Moreover, it may be possible to find a subset $\Delta_1$ of $\Omega$ which
meets each row, column and letter of $\Lambda$ in a constant number $\ell$ of
cells, where $\ell>1$, which is not a union of disjoint transversals. 
Then $\Delta_1$ is $-3$-perfect, and gives an equitable partition of $\Gamma$.
Figure~\ref{fig:4notgl} shows an example. It is not isomorphic to the one in
Figure~\ref{fig:4gl}(b), even though both have $n=4$, $k=9$, $p=3=q$ and
$a=b=6$.

\begin{figure}[htbp]
\[
\begin{array}{|c|c|c|c|}
\hline
\mathbf{A} & \mathbf{B} & C & D\\
\hline
B & C & \mathbf{D} & \mathbf{A}\\
\hline
\mathbf{C} & \mathbf{D} & A & B\\
\hline
D & A & \mathbf{B} & \mathbf{C}\\
\hline
\end{array}
\]
\caption{An equitable partition of the cyclic Latin square of order $4$, which has no 
transversal (the elements of $\Delta_1$ are shown in bold)}
\label{fig:4notgl}
\end{figure}

In a similar way we can give partitions with more than two parts, where the
parts are not unions of transversals.

\begin{eg}
Figure~\ref{fig:fano} shows the Latin square of order $7$ defined by the Steiner triple
system of order $7$.  The three different fonts show a $-3$-equitable partition with parts of
sizes $7$, $14$ and $28$.  This is strictly orthogonal to each of the
partitions into rows, columns and letters.

\begin{figure}[htbp]
\[
\begin{array}{|c|c|c|c|c|c|c|}
\hline
\mathcal{A} & \mathbf{C} & \mathbf{B} & E & D & G & F\\
\hline
C & \mathcal{B} & A & \mathbf{G} & F & E & \mathbf{D}\\
\hline
B & A & \mathcal{C} & F & \mathbf{G} & D & \mathbf{E}\\
\hline
E & G & \mathbf{F} & \mathcal{D} & A & \mathbf{C} & D\\
\hline
\mathbf{D} & F & G & \mathbf{A} & \mathcal{E} & B & C\\
\hline
G & \mathbf{E} & D & C & \mathbf{B} & \mathcal{F} & A\\
\hline
\mathbf{F} & D & E & B & C & \mathbf{A} & \mathcal{G}\\
\hline
\end{array}
\]
\caption{A Latin square of order $7$ with a $-3$-equitable partition into three parts,
one of size $7$ (calligraphic letters), one of size $14$ (bold), and one of size $28$}
\label{fig:fano}
\end{figure}
\end{eg}

\begin{eg}
If $\Delta$ is uniform (all parts have the same size)
and strictly orthogonal to each of $R$, $C$ and $L$ but
the size of the parts of $\Delta$ is not $n$, then $\{R,C,L,\Delta\}$ is a
\emph{mixed} orthogonal array.  These are discussed in \cite[Chapter 9]{oa},
whose Table 9.25 gives many examples with $n=6$ in which $\Delta$ has three
parts of size twelve.

One of these examples is shown in Figure~\ref{fig:oa}.
The natural order from \cite{oa} is used, but the underlying
Latin square is isotopic to the Cayley table of a cyclic group.

\begin{figure}[htbp]
\[
\begin{array}{|c|c|c|c|c|c|}
\hline
\mathbf{A} & \mathbf{B} & C & \mathcal{F} & \mathcal{D} & E\\
\hline
\mathcal{B} & \mathcal{C} & \mathbf{A} & D & E & \mathbf{F}\\
\hline
C & A & \mathcal{B} & \mathbf{E} & \mathbf{F} & \mathcal{D}\\
\hline
\mathbf{D} & \mathbf{E} & F & \mathcal{C} & \mathcal{A} & B\\
\hline
\mathcal{E} & \mathcal{F} & \mathbf{D} & A & B & \mathbf{C}\\
\hline
F & D & \mathcal{E} & \mathbf{B} & \mathbf{C} & \mathcal{A}\\
\hline
\end{array}
\]
\caption{A mixed orthogonal array, giving a uniform $-3$-equitable partition
  of a Latin
square of order $6$ into three parts (indicated by bold, calligraphic and normal fonts)}
\label{fig:oa}
\end{figure}

Table 12.7 of \cite{oa}
shows that Finney gave more examples  for these numbers in \cite{djf},
 and that examples with $n=10$ and parts of $\Delta$ having size $20$
were given in \cite{mandeli,wang}.
\end{eg}

\section{Mixed equitable partitions}
\label{s:mixed}

There seems even less chance of classifying equitable partitions 
of Latin-square graphs
where both
non-principal eigenvalues occur. We content ourselves with two examples.

\begin{eg}
As in any strongly regular graph, the distance partition
with respect to a vertex $\alpha$ (whose classes are $\{\alpha\}$, the 
vertices adjacent to $\alpha$, and the rest) is equitable: all three
eigenvalues occur \cite[\S4.5]{gr}.
\end{eg}

\begin{eg}
In Theorem~\ref{t:inflation_equitable}, we observed that, if $\Lambda$ is an
$s$-fold inflation of a Latin square $\Lambda_0$ of order $t$, then the partition of
$\Lambda$ into subsquares is equitable, and has $t^2$ parts. We saw that all
three eigenvalues occur if and only if $t>2$.
\end{eg}

\section{Mutually orthogonal Latin squares}
\label{s:mols}

If $\Lambda_1,\ldots,\Lambda_{m-2}$ are mutually orthogonal Latin squares of
order $n$, then we can form a graph whose vertices are the cells, two vertices
being joined if they lie in the same row or column and have the same letter in one
of the squares. This graph is strongly regular with valency $m(n-1)$
and other eigenvalues $n-m$ and $-m$. This raises the possibility of
determining the $(n-m)$-perfect sets in this graph, as a generalisation of
the main theorem of this paper. However, we expect that this will be much more
difficult.

\begin{eg}
Take two orthogonal Latin squares of order~$t=3$. The graph is complete, and
so any non-empty proper subset is $-1$-perfect, where $-1=3-4=t-m$. Hence, any
inflation of such a subset gives a $(n-4)$-perfect subset in a Latin square of
order $n=3s$ for any $s$, by the argument of
Theorem~\ref{t:inflation_equitable}.
\end{eg}

\paragraph*{Acknowledgement} The authors are grateful to Professor Yaokun Wu,
who arranged the visits of the first two authors to Shanghai Jiao Tong 
University where this research was begun. The fourth author is grateful to
Professor Denis Krotov for his interest and useful discussions.

\end{document}